\documentclass[11pt]{article}
\usepackage{flafter,amsmath,amssymb,latexsym,psfrag,graphicx,color,indentfirst}
\usepackage[margin=2.5cm]{geometry}

\usepackage[colorlinks,
linktocpage=true,
linkcolor=black,
citecolor=black,
urlcolor=black,
anchorcolor=black
]{hyperref}

\usepackage[numbers,sort&compress]{natbib}

\DeclareGraphicsExtensions{.eps,.pdf,.jpg,.png}
\allowdisplaybreaks

\newtheorem{theorem}{Theorem}[section]

\newtheorem{lemma}[theorem]{Lemma}

\numberwithin{equation}{section}
\numberwithin{figure}{section}

\begin{document}
\setlength\arraycolsep{2pt}
\date{\today}

\title{The inverse source problem for the wave equation revisited:\\ A new approach}
\author{Mourad Sini$^1$, Haibing Wang$^{2,3}$
\\$^1$RICAM, Austrian Academy of Sciences, A-4040, Linz, Austria \qquad\qquad\quad\quad
\\E-mail: mourad.sini@oeaw.ac.at
\\$^2$School of Mathematics, Southeast University, Nanjing 210096, P.R. China\;
\\$^3$Nanjing Center for Applied Mathematics, Nanjing 211135, P.R. China\qquad\,
\\E-mail: hbwang@seu.edu.cn
}

\maketitle
\begin{abstract}

The inverse problem of reconstructing a source term from boundary measurements, for the wave equation, is revisited. We propose a novel approach to recover the unknown source through measuring the wave fields after injecting small particles, enjoying a high contrast, into the medium. For this purpose, we first derive the asymptotic expansion of the wave field, based on the time-domain Lippmann-Schwinger equation. The dominant term in the asymptotic expansion is expressed as an infinite series in terms of the eigenvalues $\{\lambda_n\}_{n\in \mathbb{N}}$  of the Newtonian operator (for the pure Laplacian). Such expansions are useful under a certain scale between the size of the particles and their contrast.
Second, we observe that the relevant eigenvalues appearing in the expansion have non-zero averaged eigenfunctions. We prove that the family $\{\sin(\frac{c_1}{\sqrt{\lambda_n}}\, t),\, \cos(\frac{c_1}{\sqrt{\lambda_n}}\, t)\}$, for those relevant eigenvalues, with $c_1$ as the contrast of the small particle, defines a Riesz basis  (contrary to the family corresponding to the whole sequence of eigenvalues). Then, using the Riesz theory, we reconstruct the wave field, generated before injecting the particles, on the center of the particles. Finally, from these (internal values of these) last fields, we reconstruct the source term (by numerical differentiation for instance). A significant advantage of our approach is that we only need the measurements on $\{x\}\times(0, T)$ for a single point $x$ away from $\Omega$, i.e., the support of the source, and large enough $T$.

\bigskip
{\bf Keywords.} Inverse source problem; Asymptotic expansion; Riesz basis.\\

{\bf MSC(2010): } 35R30, 35C20, 35L05.

\end{abstract}

\section{Introduction}\label{sec_int}

Our goal in this work is to propose a different approach to study the well known inverse source problems that might have some advantages as compared to the already known ways of tackling these problems. Let us first describe the mathematical formulation of our inverse source problems.  Denote by $A_0(x)$ and $B_0(x)$  two scalar coefficients modeling the background medium. These coefficients can, respectively, model the electric permittivity $\epsilon_0$ and the magnetic permeability $\mu_0$ in the electromagnetism or the inverses of the bulk modulus $\kappa_0$ and the mass density $\rho_0$ in acoustics. Assume that they  are positive constants, i.e., $A_0(x)\equiv A_0$ and $B_0(x)\equiv B_0$. Define $c_0:=(A_0B_0)^{-1/2}$. Let $J=J(x,\,t)$ be a source compactly supported in $\Omega$, where $\Omega$ is a bounded Lipschitz domain in $\mathbb R^n$. We also assume that $J=0$ in $\Omega$ for $t<0$. Then the field $V$ generated by the source $J$ in the homogeneous background medium satisfies
\begin{equation}\label{model}
\begin{cases}
c_0^{-2} V_{tt} - \Delta V = J  &  \mathrm{in}\;\mathbb R^n\times(0,\,T),\\
V|_{t=0}=0,\;V_t|_{t=0}=0  &  \mathrm{in}\;\mathbb R^n
\end{cases}
\end{equation}
for $n=2,\, 3$.

A typical inverse source problem is to construct the source term $J$ from the measurements of $V$ on a cylinder $\partial \Omega \times (0,\,T)$ for a large time $T$. Such inverse source problems have many significant applications and have been extensively studied; see, for instance, \cite{B-Y2006, C-D1983, D-T2015, J-L-Y2017, Ton2003} and the references therein. We also refer to the monograph \cite{Isakov1990} for uniqueness and stability of various inverse source problems. There are several classical arguments to show the uniqueness and stability for inverse source problems of hyperbolic equations. The approaches based on Carleman estimates and unique continuation have been widely used \cite{B-Y2017, Klibanov1992, Tataru1996}. Based on exact boundary controllability, the stability and a reconstruction formula for the inverse problem of determining space-dependent component in the source term from boundary measurement are shown in \cite{Yama1995}.  Under a weak regularity assumption, a uniqueness result for a multidimensional hyperbolic inverse source problem with a single measurement was proved in \cite{Yama1999}. Recently, an approach of using Laplace transform was employed in \cite{Hu2020} to show the uniqueness of some inverse source problems for the wave equation.

Here, we propose a different way of looking at this problem. 
Let $D=z+a B$ be a particle injected into $\Omega$, where $a$ is a small parameter, $B\subset\mathbb R^n$ is a bounded Lipschitz domain containing the origin and $z\in\mathbb R^n$. Define
\begin{equation}\label{def_eps}
A(x):=
\begin{cases}
A_0 & \mathrm{in}\; \mathbb R^n\setminus D,\\
A_1 & \mathrm{in}\; D,
\end{cases}
\end{equation}
where $A_1$ is a positive constant with $A_1 \sim a^{-2}$ as $a\ll 1$. 
Such particles are known to exist in electromagnetism and they are called dielectric nanoparticles (with large values of the permittivity and moderate permeability). 
In acoustics, bubbles with such small values of the bulk, with moderate mass density, 
can be designed; see \cite{Z-F2018} for a discussion on this issue.

Set
$$c(x):=(A(x)B_0)^{-1/2},  \quad c_1:=(A_1B_0)^{-1/2}.$$
Define the contrast as
\begin{equation}\label{def_q}
q(x):=\frac{c_0^2}{c^2(x)}-1, \quad x\in \mathbb R^n.
\end{equation}
Then, after injecting the particle $D$ into $\Omega$, the field $U$ generated by $J$ meets
\begin{equation}\label{model_U}
\begin{cases}
c^{-2}(x) U_{tt} - \Delta U = J  &  \mathrm{in}\;\mathbb R^n\times(0,\,T),\\
U|_{t=0}=0,\;U_t|_{t=0}=0  &  \mathrm{in}\;\mathbb R^n.
\end{cases}
\end{equation}
The solvability and stability results for the direct problems \eqref{model} and \eqref{model_U} can be found in \cite[Theorem 8.1]{Isakov2017}.The aim is to reconstruct the unknown source $J$ from the measurement of $U(x,\,t)$ for a fixed $x\in\mathbb R^n\setminus \Omega$ over a large enough time interval, that is, measuring the wave fields after injecting the small particle $D$. The results will be given and proved in the dimensions $n=3$. The corresponding results and proofs for the case $n=2$ can be done similarly at the expense of handling the different type of the singularity of the Green's function.

We first analyze the asymptotic behavior of the wave field $U(x,\,t)$ for the wave equation model as $a\to 0$. Asymptotic expansions of physical fields generated by small particles are well developed in the literature for elliptic models, but there are few results on time-domain models. Recently, we studied a heat conduction problem and an acoustic scattering problem with many small holes in \cite{S-W2019} and \cite{S-W-Y2021}, respectively, where we derived the asymptotic analysis of the solutions  by using time-domain integral equation methods. In this paper, based on the time-domain Lippmann-Schwinger equation, we shall show the asymptotic expansion of $U(x,\,t)$ for $x\in\mathbb R^3\setminus\overline D$ as $a\to 0$ through the eigensystem of the Newtonian potential operator related to the Laplace equation. We set $W:=U-V$. Then we have
\begin{equation}\label{model_W}
\begin{cases}
c_0^{-2} W_{tt} - \Delta W = -\dfrac{q}{c_0^2} U_{tt}  &  \mathrm{in}\;\mathbb R^3\times(0,\,T),\\
W|_{t=0}=0,\;W_t|_{t=0}=0  &  \mathrm{in}\;\mathbb R^3.
\end{cases}
\end{equation}
The result is as follows.

\begin{theorem}\label{Asymptot-Introduction} Under the conditions above, we have the expansion
\begin{eqnarray*}
W(x,\,t) &=&\sum_{n=1}^{+\infty}\frac{c_1\lambda_n^{-\frac{3}{2}}}{4\pi |x-z|}\,\left(\int_D e_n(y)\,dy\right)^2\,\int_{c_0^{-1}|x-z|}^{t} \sin\left[ \frac{c_1}{\sqrt{\lambda_n}}(t-s)\right]\,V(z,\,s-c_0^{-1}|x-z|)\,ds \nonumber\\
&& \qquad\qquad\qquad\qquad\qquad\qquad + O(a^2), \quad (x,\,t) \in (\mathbb R^3\setminus\overline D)\times (0,\,T),
\end{eqnarray*}
where 
$\{\lambda_n, e_n\}_{n\in \mathbb{N}}$ is the family of the eigenelements of the Newtonian operator $N:\, L^2(D)\rightarrow L^2(D)$ defined by 
\begin{equation}\label{op_N}
N(f)(x):=\int_{D}\frac{f(y)}{4 \pi\vert x-y\vert}d(y).
\end{equation}
\end{theorem}

Observe that the relevant eigenfunctions $e_n$'s are those with non-vanishing averages. Note that for the Newtonian operator, $1$ is not an eigenvalue. Therefore, the family of eigenelements $\{\lambda_n,\, e_n\}$, with non-vanishing averages, is rich. In addition, in the case that $D$ is a ball, the behavior in terms of the index $n$, i.e., the Weyl law, of these relevant eigenvalues, that we still denote $\{\lambda_n\}_{n \in \mathbb{N}}$, allow us to show that the family $\{\sin(\frac{c_1}{\sqrt{\lambda_n}}\, t),\,\cos(\frac{c_1}{\sqrt{\lambda_n}}\, t)\}$, for those relevant eigenvalues, defines a Riesz basis (contrary to the family corresponding to the whole sequence of eigenvalues). We state this as a second result.

\begin{theorem}\label{Riesz-Introduction} 
Let $D$ be a ball and the set $\{\lambda_n, e_n\}_{n\in \mathbb{N}}$ be the subfamily of the eigenelements of the Newtonian operator $N$ satisfying $\int_B e_n(x)dx\neq 0$. Then the family $\{\sin(\omega_n\, t), \, \cos(\omega_n\, t)\}_{n\in\mathbb N}$
defines a Riesz basis in $L^2[-\frac{\pi}{b},\, \frac{\pi}{b}]$, where $b=\frac{c_1\pi}{a}$ and $\omega_n:=\frac{c_1}{\sqrt{\lambda_n}},\; n \in \mathbb{N}$.
\end{theorem}

We do believe that this result can be extended to smooth general shapes.
\bigskip

The dominant term in the expansion is expressed as an infinite series in terms of Riesz basis. Using the Riesz theory, we can approximately obtain the Riesz coefficients of the wave field $V$ at $z$ over a time interval. This is the object of the final result in this introduction.

\begin{theorem}\label{Reconstruction-Introduction}
Assume that the source $J$ is compactly supported in $\Omega\times[0,\,T]$. Let $x\in \mathbb R^3\setminus \Omega$ and $z \in \Omega$, i.e., the center of the particle, be fixed. We measure both $V(x,\, t)$ and $U(x,\,t)$ for an interval of the form $[\tilde{T},\, \tilde{T}+\frac{2\pi}{b}]$ with $\tilde{T}>T_J+\frac{\vert x-z \vert}{c_0}$ and $T_J:=T+\max_{y\in \Omega}\frac{\vert z-y \vert}{c_0}$. Then
\begin{equation}\label{Riesz_h7}
V(z,\,t) = \sum_{n=1}^{+\infty}\left[C_n(x,\,z)g_n(t-\pi/b) + D_n(x,\,z) h_n(t-\pi/b)\right], \quad t\in[0,\,2\pi/b]
\end{equation}
with \begin{eqnarray}
C_n(x,\,z)&:=& \int_{0}^{2\pi/b} \cos\left(\omega_n (s-\pi/b)\right) \,V(z,\,s)\,ds,\quad n\in\mathbb N, \label{Fourier-modes-1}\\
D_n(x,\,z)&:=& \int_{0}^{2\pi/b} \sin\left(\omega_n (s-\pi/b)\right) \,V(z,\,s)\,ds,\quad n\in\mathbb N, \label{Fourier-modes-2}
\end{eqnarray}
where $g_n(t):=(\mathcal T^{-1})^*\cos((n-1) t)$ and  $h_n(t):=(\mathcal T^{-1})^*\sin(n t)$.
The operator $\mathcal T$ is the one linking the Riesz basis $\{\sin(\omega_n\, t), \,\cos(\omega_n\, t)\}_{n\in\mathbb N}$ to the canonical basis $\{\sin(n\, t), \,\cos((n-1)\, t)\}_{n\in\mathbb N}$; see (\ref{defn_T}).
The coefficients $C_n(x, \,z)$ and $D_n(x, \,z)$ can be estimated from the data $W(x,\, t)$ as follows:
\begin{eqnarray*}
C_n(x,\,z) &=& \cos(\omega_n\tilde T) A_n(x,\,z) -  \sin(\omega_n\tilde T) B_n(x,\,z),\\
D_n(x,\,z) &=& \sin(\omega_n\tilde T) A_n(x,\,z) + \cos(\omega_n\tilde T) B_n(x,\,z),
\end{eqnarray*}
where
\begin{eqnarray}
A_n(x,\,z)&=& \alpha^{-1}_n \int_{0}^{2\pi/b} h_n(t-\pi/b) \, U(x,\,\tilde T+s)\,ds +O(a),\quad n\in\mathbb N, \label{Fourier-modes-1-computed}\\
B_n(x,\,z)&=& \alpha^{-1}_n \int_{0}^{2\pi/b} g_n(t-\pi/b) \, U(x,\,\tilde T+s)\,ds +O(a),\quad n\in\mathbb N \label{Fourier-modes-2-computed}
\end{eqnarray}
with
\begin{equation}\label{alpha-1}
\alpha_n:=\alpha_n(x,\,z):= \frac{\omega_n \left(\int_D e_n(y)\,dy\right)^2}{4\pi |x-z| \lambda_n}.
\end{equation}
\end{theorem}

Varying the location of $z$ in $\Omega$, we get the wave field $V$ in $\Omega$ over a time interval. Then, we can further reconstruct the source term $J$ by numerically differentiating $V$. A significant advantage of the proposed method is that we only need to perform the measurement of $U$ at a fixed point $x\in\mathbb R^3\setminus\Omega$, which is attractive in practical situation. As we reconstruct $V(z,\, t)$ for $ t\in[0,\,2\pi/b]$, then if $2\pi/b > T_J$, we can reconstruct $V$ for any time. This condition on $b=c_1\pi/a$ makes sense since it is related to the used particle which is at our handle.

This approach is inspired by the imaging modalities using contrast agents; see \cite{B-B-I2011, S-K-V-H2010, McLeod-Ozcan, Ntziachristos} and the references therein for real world applications. The common motivation there is that the imaging techniques are potentially able to extract features when the physical parameters between healthy tissues and malignant ones have relatively high contrast. However, it is observed that for a benign tissue, the variation of the permittivity is quite low so that conventional imaging modalities are limited to be employed for early detection of diseases. In such cases, it is highly desirable to create such missing contrast. One way to do is to use electromagnetic nanoparticles as contrast agents. The idea can also be used to acoustic imaging modalities using bubbles as contrast agents. The mathematical analysis of such imaging modalities are given in \cite{C-C-S2018, D-G-S2021, G-S2022, G-S2022-2}. The main idea there is that contrasting the fields measured before and after injecting the contrasting particles, one can reconstruct the total fields on the location of those particles. In \cite{D-G-S2021, G-S2022, G-S2022-2}, this is justified in the time harmonic regimes (both for the electromagnetic and acoustic models)  when the used incident waves are sent at frequencies close to the resonating frequencies of the used particles. These resonances exist for certain scales of the contrasts of the particles. In our work here, we extend such ideas to the time domain inverse source problem. In this time-domain regime, the same scales of the contrasts allow us to extract the Fourier modes, in time, of the internal field using the eigenvalues of the Newtonian operator; see (\ref{Fourier-modes-1})--(\ref{Fourier-modes-2}). Observe that for these scales the resonances, called the dielectric resonances, used in \cite{D-G-S2021, G-S2022} are related, and actually due, to the eigenvalues of the Newtonian operator as well.

We finish this introduction by the following observation. The expansion in Theorem \ref{Asymptot-Introduction} means that the field $U(x, t)$ is composed of the field $V(x, t)$, which we call the first field, generated by solely the source $f$ and the one generated by the interaction of this field and particle $D$ which we call the second field. This second field, which is a local spot, is modeled by (\ref{L-S_3}); see also (\ref{L-S_6}). As mentioned in Theorem \ref{Reconstruction-Introduction}, while measuring $U(x, t)$ for large $t$, the first wave field $V(x, t)$ has already fully passed the point $x$ where we measure. Therefore, we only measure the second field due to the interaction between the source $f$ and the small particle. 

The remainder of the paper is organized as follows. In Section \ref{sec_asym}, we derive the time-domain Lippmann-Schwinger equation for $W$ and then prove Theorem \ref{Asymptot-Introduction}. In Section \ref{sec-Riesz}, we prove Theorem \ref{Riesz-Introduction} and finally in Section \ref{sec-Reconstruction}, we prove Theorem \ref{Reconstruction-Introduction}.

\section{Proof of Theorem \ref{Asymptot-Introduction}}\label{sec_asym}

In this section, we show the asymptotic behavior of the field $W(x,\,t)$ for $(x,\,t)\in (\mathbb R^3\setminus\overline D)_T$ as $a\ll 1$, and completes the proof of Theorem  \ref{Asymptot-Introduction}.

To proceed, we introduce the function space
\begin{equation*}
H_0^r(0,\,T):=\left\{ g|_{(0,\,T)}:\; g\in H^r(\mathbb R)\, \textrm{ with }\, g\equiv 0 \, \textrm{ in }\, (-\infty,\,0) \right\},\quad r\in\mathbb R,
\end{equation*}
and generalize it to the $E$-valued function space, denoted by $H_0^r(0,\,T;\, E)$, where $E$ is a Hilbert space. Define
\begin{equation*}
LT(\sigma,\,E):= \left\{f\in \mathcal D^\prime_+(E):\,e^{-\sigma t}f \in \mathcal S^\prime_+(E) \right\}, \quad \sigma>0,
\end{equation*}
where $\mathcal D^\prime_+(E)$ and $\mathcal S^\prime_+(E)$ denote the sets of distributions and temperate distributions on $\mathbb R$, respectively, with values in $E$ and supports in $[0,\,+\infty)$. Then we define the space
\begin{equation*}
H_{0,\sigma}^r(0,\,T;\, E):= \left\{f\in LT(\sigma,\,E) : \, e^{-\sigma t}\partial_t^r f\in L^2_0(0,\,T;\,E)  \right\}, \quad r\in\mathbb Z_+.
\end{equation*}
For nonnegative integer $r$, we introduce the norm
\begin{equation}\label{norm}
\| f \|_{H_{0,\sigma}^r(0,\,T;\, E)}:=\left(\int_0^T  e^{-2\sigma t}\left[\left\| f \right\|^2_E + \sum_{k=1}^r T^{2k} \left\| \frac{\partial^k f}{\partial t^k} \right\|^2_E\right] \,dt\right)^{1/2}.
\end{equation}
For simplicity of notations, we denote $X\times (0,\,T)$ by $X_T$, where $X$ is a domain in $\mathbb R^3$. We shall also use the notation \lq\lq$\lesssim$\rq\rq{} to denote \lq\lq$\leq$\rq\rq{} with its right-hand side multiplied by a generic positive constant, if we do not emphasize the dependence of the constant on some parameters.

For convenience, we set
\begin{equation*}
q_0=\frac{c_0^2}{c_1^2}-1 \sim a^{-2} \quad \textrm{for}\;a\ll 1.
\end{equation*}
Define the retarded volume potential $V_D$ by
\begin{equation}\label{op_V_D}
V_D[f](x,\,t):=\int_D \frac{f(y,\,t-c_0^{-1}|x-y|)}{4\pi|x-y|}\,dy, \quad (x,\,t)\in \mathbb R^3\times(0,\,T).
\end{equation}
From \eqref{model_W}, we are led to the Lippmann-Schwinger equation
\begin{equation}\label{L-S}
U(x,\,t) + \frac{q_0}{c_0^2} V_D[U_{tt}](x,\,t)= V(x,\,t), \quad (x,\,t)\in \mathbb R^3\times(0,\,T)
\end{equation}
or
\begin{equation}\label{L-S_0}
W(x,\,t) + \frac{q_0}{c_0^2}V_D[W_{tt}](x,\,t) = -\frac{q_0}{c_0^2}V_D[V_{tt}](x,\,t), \quad (x,\,t)\in \mathbb R^3\times(0,\,T).
\end{equation}
Following the convolution quadrature based argument in \cite{Lubich1994}, we can prove that the equation \eqref{L-S_0} has a unique solution in $H_{0,\sigma}^r(0,\,T;\,L^2(D))$ for $V\in H_{0,\sigma}^{r+2}(0,\,T;\,L^2(D))$.
%Moreover, we have the following {\it a-priori} estimate
%\begin{equation}\label{est_sigma_1}
%\| u^s \|_{H_{0,\sigma}^r(0,\,T;\,L^2(D))} \lesssim \| V \|_{H_{0,\sigma}^{r+2}(0,\,T;\,L^2(D))},\quad r\in\mathbb R.
%\end{equation}
%By the embedding $H^r(0,\,T)\hookrightarrow C[0,\,T]$ for $r>1/2$, we also have
%\begin{equation}\label{est_sigma_2}
%\| u^s(\cdot,\,t) \|_{L^2(D)} \lesssim \| V \|_{H_{0,\sigma}^r(0,\,T;\,L^2(D))}, \quad r>5/2,\;t\in [0,\,T].
%\end{equation}

\bigskip
In the following, we derive the asymptotic behavior of the field $W(x,\,t)$. We start by scaling both the space and time variables. Set
$$T_a := T/a.$$
For any functions $\varphi$ and $\psi$ defined on $D_T$ and $B_{T_a }$, respectively, we use the notations
\begin{eqnarray*}
&&\hat\varphi(\xi,\,\tau)=\varphi^\wedge(\xi,\,\tau):=\varphi(a \xi+z,\,a  \tau), \quad (\xi,\,\tau)\in B_{T_a },\\
&&\check\psi(x,\,t)=\psi^\vee(x,\,t):=\psi\left(\frac{x-z}{a },\,\frac{t}{a }\right), \quad (x,\,t)\in D_T.
\end{eqnarray*}
Notice that
\begin{equation}\label{scaling_t}
\frac{\partial^n \hat\varphi(\cdot,\,\tau)}{\partial \tau^n}=\frac{\partial^n \varphi(\cdot,\,a  \tau)}{\partial \tau^n} = a ^n \frac{\partial^n \varphi(\cdot,\,t)}{\partial t^n}, \quad n \in\mathbb Z_+.
\end{equation}
Then we have the following scaling result.

\begin{lemma}\label{scalingt}
Suppose $0< a  \leq 1$. For $\psi \in H^r_{0,\sigma}(0,\,T;\, L^2(D))$  with nonnegative integer $r$, we have
\begin{equation}\label{scalingt_11}
\|\psi\|_{H^r_{0,\sigma}(0,\,T;\, L^2(D))} = a^2 \|\hat \psi\|_{H^r_{0,a \sigma}(0,\,T_a ;\, L^2(B))}.
\end{equation}
\end{lemma}

{\bf Proof.} From \eqref{norm} and \eqref{scaling_t}, we can easily derive that
\begin{eqnarray*}
\|\psi\|^2_{H^r_{0,\sigma}(0,\,T;\,L^2(D))} &=&  \int_0^T e^{-2\sigma t}\, \sum_{k=0}^r T^{2k} \,\left\|\frac{\partial^k \psi(\cdot,\,t)}{\partial t^k}\right\|^2_{L^2(D)}\,dt\\
& = & a^4 \int_0^{T_a } e^{-2\sigma a  \tau}\, \sum_{k=0}^r T_a^{2k} \,\left\|\frac{\partial^k \hat\psi(\cdot,\,\tau)}{\partial \tau^k}\right\|^2_{L^2(B)}\, d\tau \\
& = & a ^4 \, \|\hat \psi\|^2_{H^r_{0,a \sigma}(0,\,T_a ;\,L^2(B))}.
\end{eqnarray*}
The proof is complete. \hfill $\Box$

\bigskip
Let
\begin{equation}\label{op_S_D}
\mathcal S_D:=-(I + (q_0/c_0^2)\,V_D\,\partial_t^2)^{-1}\, (q_0/c_0^2)\,V_D\,\partial_t^2.
\end{equation}
Define the retarded volume potential $V_B^a$ by
\begin{equation}\label{op_V_B}
V_B^a[\psi](\xi,\,\tau):=\int_B \frac{\psi(\eta,\,\tau-c_0^{-1}|\xi-\eta|)}{4\pi|\xi-\eta|}\,d\eta, \quad (\xi,\,\tau)\in \mathbb R^3\times(0,\,T_a),
\end{equation}
and then introduce the operator $\mathcal S_B^a$ by
\begin{equation}\label{op_S_B}
\mathcal S_B^a:=-(I + (q_0/c_0^2)\,V_B^a\,\partial_\tau^2)^{-1}\, (q_0/c_0^2)\,V_B^a\,\partial_\tau^2.
\end{equation}
Then we can easily prove the following result.
\begin{lemma}\label{scalingV}
Let $\psi\in H^{r+2}_{0,\sigma}(0,\,T;\,L^2(D))$. Then
\begin{eqnarray}
V_D[\psi_{tt}] &=& (V^a_B[\hat\psi_{\tau\tau}])^\vee,  \label{scalingV_1}\\
\|\mathcal S_D\|_{\mathcal L\left(H^{r+2}_{0, \sigma}(0,\,T;\,L^2(D)),\, H^r_{0, \sigma}(0,\,T;\,L^2(D))\right)} &=& \|\mathcal S_B^a\|_{\mathcal L\left(H^{r+2}_{0,a \sigma}(0,\,T_a ;\,L^2(B)),\, H^r_{0,a \sigma}(0,\,T_a ;\,L^2(B))\right)}.   \label{scalingS_1}
\end{eqnarray}
\end{lemma}

{\bf Proof.} Let $x=a \xi + z$, $y=a \eta + z$ and $t=a \tau$. Then we have
\begin{eqnarray*}
V_D[\psi_{tt}](x,\,t) &=& \int_D \frac{\psi_{tt}(y,\,t-c_0^{-1}|x-y|)}{4\pi |x-y|}\,dy\\
&=& \int_B \frac{\psi_{tt}(a \eta+z,\,a \tau-c_0^{-1}a|\xi-\eta|)}{4\pi a |\xi-\eta|}\,a^3\,d\eta\\
&=& \int_B \frac{a^{-2}\hat\psi_{\tau\tau}(\eta,\, \tau-c_0^{-1}|\xi-\eta|)}{4\pi a |\xi-\eta|}\,a^3\,d\eta\\
&=& V^a_B[\hat \psi_{\tau\tau}](\xi,\,\tau),
\end{eqnarray*}
which completes the proof of \eqref{scalingV_1}.

Further, we have $(\mathcal S_D[\psi])^\wedge=\mathcal S_B^a[\hat \psi]$. Then, it can be derived that
\begin{eqnarray*}
& &\left \| \mathcal S_D \right\|_{\mathcal L \left(H^{r+2}_{0,\sigma}(0,\,T;\,L^2(D)),\,H^r_{0,\sigma}(0,\,T;\,L^2(D))\right)}\\
&:=& \sup\limits_{0\not=\psi\in H^{r+2}_{0,\sigma}(0,\,T;\,L^2(D))} \displaystyle \frac{ \| \mathcal S_D[\psi] \|_{H^r_{0,\sigma}(0,\,T;\,L^2(D))}}{\|\psi\|_{H^{r+2}_{0,\sigma}(0,\,T;\,L^2(D))}}\\
&=& \sup\limits_{0\not=\psi\in H^{r+2}_{0,\sigma}(0,\,T;\,L^2(D))} \displaystyle \frac{ a^2\, \| \left(\mathcal S_D[\psi]\right)^\wedge \|_{H^r_{0,a\sigma}(0,\,T_a;\,L^2(B))}}{a^2 \|\hat\psi\|_{H^{r+2}_{0,a\sigma}(0,\,T_a;\,L^2(B))}}\\
&=& \sup\limits_{0\not=\hat\psi\in H^{r+2}_{0,a\sigma}(0,\,T_a;\,L^2(B))}  \displaystyle \frac{ \| \mathcal S^a_B[\hat\psi] \|_{H^r_{0,a\sigma}(0,\,T_a;\,L^2(B))}}{\|\hat\psi\|_{H^{r+2}_{0,a\sigma}(0,\,T_a;\,L^2(B))}}\\
&=& \left \| \mathcal S^a_B \right\|_{\mathcal L \left(H^{r+2}_{0,a\sigma}(0,\,T_a;\,L^2(B)),\,H^r_{0,a\sigma}(0,\,T_a;\,L^2(B))\right)}.
\end{eqnarray*}
The proof is complete. \hfill $\Box$

\medskip
Denote by $\hat V_B$ the volume potential operator for the Helmholtz equation
$$-\Delta U + s^2\,U=\hat f.$$
Define $A(s):\,L^2(B)\to L^2(B)$ by
\begin{equation}\label{def_A}
A(s):=-(I + (q_0s^2/c_0^2)\,\hat V_B)^{-1}\, (q_0s^2/c_0^2)\,\hat V_B,
\end{equation}
which is the Fourier-Laplace transform of $\mathcal S_B^a$ defined by \eqref{op_S_B}. Then we have the following estimate in \cite[Theorem 4.1]{L-M2015}:
\begin{equation}\label{esi_A}
\left\| A(s) \right\|_{\mathcal L(L^2(B),\,L^2(B))} \leq  \frac{q_0\,|s|^2}{\sigma^2\,c_0^2}
\end{equation}
for any $s$ with $\mathrm{Re}\, s = \sigma>0$.

\begin{lemma}\label{inv_SB}
The operator $\mathcal S_B^a :\, H^{r+2}_{0,a \sigma}(0,\,T_a ;\,L^2(B)) \to H^r_{0,a \sigma}(0,\,T_a ;\,L^2(B))$ with $r=0,\,1$ is estimated by $O(q_0)=O(a^{-2})$ for $a \ll 1$.
\end{lemma}

{\bf Proof.} Note that $\mathcal S_B^a$ is the corresponding time-domain solution operator obtained by the inverse Fourier-Laplace transform of $A(s)$. Using Lubich's notation, we have $\mathcal S_B^a = A(\partial_t)$. Then, for $r=0$ and $g\in H^2_{0,a \sigma}(0,\,T_a ;\,L^2(B))$ compactly supported with respect to $t$ in $[0,\,T_a ]$, we have
\begin{eqnarray*}\label{SB_2}
\left\| \mathcal S^a_B[g]\right\|^2_{H^0_{0,a \sigma}(0,\,T_a ;\,L^2(B))} &=&\int_0^{T_a } \left(e^{-a \sigma t} \,\|A(\partial_t)g \|_{L^2(B)}\right)^2\,dt\nonumber \\
& \leq & \int_0^{+\infty}\left(\mathcal F \left[e^{-a \sigma t} \|A(\partial_t)g \|_{L^2(B)}\right](\eta)\right)^2\,d\eta \nonumber\\
& = & \int_{a \sigma-i\mathbb R_+}\left(\mathcal L \left[ \|A(\partial_t)g \|_{L^2(B)}\right](s)\right)^2\,ds \nonumber\\
& = & \int_{a \sigma-i\mathbb R_+}\|\mathcal L \left[ A(\partial_t)g\right](s) \|^2_{L^2(B)}\,ds \nonumber\\
& = & \int_{a \sigma-i\mathbb R_+}\| A(s)\,\mathcal L[g](s) \|^2_{L^2(B)}\,ds \nonumber\\
& \leq & \frac{q_0^2}{c_0^4\,(a \sigma)^4}\,\int_{a \sigma-i\mathbb R_+}\| s^2 \mathcal L[g](s) \|^2_{L^2(B)}\,ds \nonumber\\
& = & \frac{q_0^2}{c_0^4\,(a \sigma)^4}\,\int_{a \sigma-i\mathbb R_+} \| \mathcal L[\partial_t^2 g](s) \|^2_{L^2(B)}\,ds \nonumber\\
& = & \frac{q_0^2\,T_a^{-4}}{c_0^4\,(a \sigma)^4}\,\int_0^{+\infty} \left( e^{-a \sigma t}\,T_a^2\,\| \partial_t^2 g \|_{L^2(B)}\right)^2\,dt \nonumber\\
& \leq & \frac{q_0^2\,T^{-4}}{c_0^4\,\sigma^4}\,\left\| g\right\|^2_{H^2_{0,a \sigma}(0,\,T_a ;\,L^2(B))} .
\end{eqnarray*}
Similarly, we derive for $r=1$ that
\begin{eqnarray*}\label{SB_3}
\left\| \mathcal S^a_B[g]\right\|^2_{H^1_{0,a \sigma}(0,\,T_a ;\,L^2(B))}&=& \int_0^{T_a }e^{-2a \sigma t} \left( \| A(\partial_t)g \|^2_{L^2(B)} + T_a^2\, \left\|\frac{\partial A(\partial_t)g}{\partial t}\right\|^2_{L^2(B)}\right)\,dt \nonumber\\
&\leq & 2T_a ^2\, \int_0^{T_a }e^{-2a \sigma t}\, \|\partial_t\left( A(\partial_t)g\right)\|^2_{L^2(B)}  \,dt\nonumber\\
& \leq & 2T_a ^2\, \int_0^{+\infty} \left(\mathcal F \left[e^{-a \sigma t}\, \| \partial_t\left(A(\partial_t)g\right)\|_{L^2(B)} \right](\eta)\right)^2 \,d\eta \nonumber\\
& \leq &  \frac{2T_a ^2\,q_0^2}{c_0^4\,(a \sigma)^4}\,\int_{a \sigma-i\mathbb R_+} \|s^3 \mathcal L[g](s) \|^2_{L^2(B)}\,ds \nonumber\\
& = & \frac{2T_a ^2\,q_0^2}{c_0^4\,(a \sigma)^4}\,\int_{a \sigma-i\mathbb R_+}\| \mathcal L[\partial_t^3 g](s) \|^2_{L^2(B)}\,ds \nonumber\\
& = & \frac{2T_a^{-4}\,q_0^2}{c_0^4\,(a \sigma)^4}\,\int_0^{+\infty} e^{-2a \sigma t}\, T_a^6\,\| \partial_t^3 g \|^2_{L^2(B)}\,dt\nonumber \\
&\leq & \frac{2T^{-4}\,q_0^2}{c_0^4\,\sigma^4}\,\left\|g\right\|^2_{H^3_{0,a \sigma}(0,\,T_a ;\,L^2(B))}.
\end{eqnarray*}
The proof is complete. \hfill $\Box$

\bigskip
Now we estimate $W|_{D_T}$ by
\begin{eqnarray}\label{est_us_D}
\|W(\cdot,\,t)\|_{L^2(D)} &\lesssim& \|W\|_{H^1_{0,a \sigma}(0,\,T;\,L^2(D))} \nonumber\\
&\lesssim& \|\mathcal S_D\|_{\mathcal L\left(H^3_{0, \sigma}(0,\,T;\,L^2(D)),\, H^1_{0, \sigma}(0,\,T;\,L^2(D))\right)}\, \|V\|_{H^3_{0, \sigma}(0,\,T;\,L^2(D))} \nonumber\\
&\lesssim& a^{3/2}\, \|\mathcal S_B\|_{\mathcal L\left(H^3_{0,a \sigma}(0,\,T_a ;\,L^2(B)),\, H^1_{0,a \sigma}(0,\,T_a ;\,L^2(B))\right)} \nonumber\\
&\lesssim& a^{-1/2},
\end{eqnarray}
and hence
\begin{equation}\label{est_u_D}
\|U(\cdot,\,t)\|_{L^2(D)} \lesssim a^{-1/2}
\end{equation}
for any $t\in[0,\,T]$. In the same way, we can also estimate $\partial_t^k U|_{D_T}$ by
\begin{equation}\label{est_u_D_t}
\|\partial_t^k U(\cdot,\,t)\|_{L^2(D)} \lesssim a^{-1/2}, \quad t\in[0,\,T],\;k=1,\,2,\,\cdots.
\end{equation}

We rewrite \eqref{L-S} as
\begin{eqnarray}\label{L-S_2}
U(x,\,t) + \frac{q_0}{c_0^2}\int_D \frac{U_{tt}(y,\,t)}{4\pi|x-y|}\,dy &=& V(x,\,t) + \frac{q_0}{c_0^2}\int_D \frac{U_{tt}(y,\,t)-U_{tt}(y,\,t-c_0^{-1}|x-y|)}{4\pi|x-y|}\,dy \nonumber\\
&=:& V(x,\,t) + O(q_0a), \quad (x,\,t)\in \mathbb R^3\times(0,\,T).
\end{eqnarray}
Consider the following problem:
\begin{equation}\label{L-S_3}
\begin{cases}
\varphi(x,\,t) + \displaystyle\frac{q_0}{c_0^2}\int_D \frac{\varphi_{tt}(y,\,t)}{4\pi|x-y|}\,dy = f(x,\,t) & \mathrm{in}\;D\times(0,\,T),\\
\varphi(x,\,0)=\varphi_t(x,\,0)=0 & \mathrm{in}\;D.
\end{cases}
\end{equation}
Let $\{\lambda_n,\,e_n\}$ be the eigensystem of the Newtonian potential operator $N$ defined by
\begin{equation*}
N[\varphi](x,\,t):=\int_D \frac{\varphi(y,\,t)}{4\pi |x-y|}\,dy, \quad (x,\,t)\in D_T.
\end{equation*}
Set
\begin{equation*}
\varphi_n(t):=\langle \varphi(\cdot,\,t),\, e_n\rangle, \quad f_n(t):=\langle f(\cdot,\,t),\, e_n\rangle.
\end{equation*}
We have
\begin{equation*}
\varphi(x,\,t) = \sum_{n=1}^{+\infty} \langle \varphi(\cdot,\,t),\, e_n\rangle\, e_n = \sum_{n=1}^{+\infty} \varphi_n(t)\, e_n.
\end{equation*}
Then it is derived from \eqref{L-S_3} that
\begin{equation}\label{L-S_4}
\begin{cases}
\frac{q_0}{c_0^2}\lambda_n \varphi_n^{\prime\prime}(t) + \varphi_n(t) = f_n(t) \quad \mathrm{in}\; (0,\,T),\\
\varphi_n(0)=\varphi_n^\prime(0)=0
\end{cases}
\end{equation}
for $n=1,\,2,\,\cdots$. The solution to \eqref{L-S_4} is given by
\begin{equation}\label{L-S_5}
\varphi_n(t)=c_0 (q_0\lambda_n)^{-1/2} \int_0^t \sin\left[c_0(q_0\lambda_n)^{-1/2}(t-\tau)\right]\,f_n(\tau)\,d\tau.
\end{equation}
Moreover, we have the estimate
\begin{equation}\label{est_phi}
\|\varphi(\cdot,\,t)\|_{L^2(D)}\lesssim \|f(\cdot,\,t)\|_{L^2(D)}.
\end{equation}
Then we get from \eqref{L-S_2} that
\begin{equation}\label{est_u_D2}
\|U(\cdot,\,t)\|_{L^2(D)} \lesssim a^{1/2}
\end{equation}
for any $t\in[0,\,T]$. In the same way, we can also improve the estimate of \eqref{est_u_D_t} as
\begin{equation}\label{est_u_D_t2}
\|\partial_t^k U(\cdot,\,t)\|_{L^2(D)} \lesssim a^{1/2}, \quad t\in[0,\,T],\;k=1,\,2,\,\cdots.
\end{equation}
Then the error term in \eqref{L-S_2} can be improved as
\begin{eqnarray}\label{L-S_22}
U(x,\,t) + \frac{q_0}{c_0^2}\int_D \frac{U_{tt}(y,\,t)}{4\pi|x-y|}\,dy &=& V(x,\,t) + \frac{q_0}{c_0^2}\int_D \frac{U_{tt}(y,\,t)-U_{tt}(y,\,t-c_0^{-1}|x-y|)}{4\pi|x-y|}\,dy \nonumber\\
&=:& V(x,\,t) + O(q_0a^2) , \quad (x,\,t)\in \mathbb R^3\times(0,\,T).
\end{eqnarray}
Using the estimate \eqref{est_phi} for \eqref{L-S_22}, we get
\begin{equation}\label{est_u_D3}
\|U(\cdot,\,t)\|_{L^2(D)} \lesssim a^{3/2}.
\end{equation}
Similarly, we can have
\begin{equation}\label{est_u_D_t3}
\|\partial_t^k U(\cdot,\,t)\|_{L^2(D)} \lesssim a^{3/2}, \quad t\in[0,\,T],\;k=1,\,2,\,\cdots.
\end{equation}

Let $\tilde U$ be the solution to \eqref{L-S_3} with
$f(x,\,t)=V(z,\,t)$. Then we have
\begin{equation*}
f_n(t)=V(z,\,t)\int_D e_n(x)\,dx,
\end{equation*}
and hence
\begin{equation}\label{L-S_6}
\tilde U(x,\,t)=\sum_{n=1}^{+\infty}c_0 (q_0\lambda_n)^{-1/2}\, e_n(x)\, \int_0^t \sin\left[c_0(q_0\lambda_n)^{-1/2}(t-\tau)\right]\,V(z,\,\tau)\,d\tau\,\int_D e_n(y)\,dy
\end{equation}
for $(x,\,t)\in D\times (0,\,T)$. Thus, we have
\begin{eqnarray}\label{L-S_7}
U(x,\,t)&=&\tilde U(x,\,t) + O(a^2) \nonumber\\
&=&\sum_{n=1}^{+\infty}c_0 (q_0\lambda_n)^{-1/2}\, e_n(x)\, \int_0^t \sin\left[c_0(q_0\lambda_n)^{-1/2}(t-\tau)\right]\,V(z,\,\tau)\,d\tau\,\int_D e_n(y)\,dy \nonumber\\
& & \qquad + O(a^2).
\end{eqnarray}

For $x$ away from $D$ and $y\in D$, since
\begin{equation*}
\left|\frac{1}{|x-z|} - \frac{1}{|x-y|}\right|=O(a)
\end{equation*}
and
\begin{equation*}
\left|U_{tt}(y,\,t-c_0^{-1}|x-z|) - U_{tt}(y,\,t-c_0^{-1}|x-y|)\right|=\left|\partial_t^3 U(y,\,t^*) \right|\, |y-z|,
\end{equation*}
we have
\begin{eqnarray}\label{L-S_1}
W(x,\,t)&=& -\frac{q_0}{4\pi c_0^2 |x-z|}\int_D U_{tt}(y,\,t-c_0^{-1}|x-z|)\,dy \nonumber\\
& & + \frac{q_0}{4\pi c_0^2}\int_D \left( \frac{1}{|x-z|} - \frac{1}{|x-y|} \right)\,U_{tt}(y,\,t-c_0^{-1}|x-z|)\,dy \nonumber\\
& & + \frac{q_0}{c_0^2}\int_D \frac{U_{tt}(y,\,t-c_0^{-1}|x-z|) - U_{tt}(y,\,t-c_0^{-1}|x-y|)}{4\pi|x-y|} \,dy \nonumber\\
&=& -\frac{q_0}{4\pi c_0^2 |x-z|}\int_D U_{tt}(y,\,t-c_0^{-1}|x-z|)\,dy + O(q_0a)\int_D U_{tt}(y,\,t-c_0^{-1}|x-z|)\,dy\nonumber\\
&&+ O(q_0a)\int_D \left|\partial_t^3 U(y,\,t^*) \right|\,dy \nonumber\\
&=& -\frac{q_0}{4\pi c_0^2 |x-z|}\int_D U_{tt}(y,\,t-c_0^{-1}|x-z|)\,dy + O(a^2).
\end{eqnarray}
Inserting \eqref{L-S_7} into \eqref{L-S_1}, we obtain
\begin{eqnarray}\label{L-S_8}
&&W(x,\,t) \nonumber\\
&=&\sum_{n=1}^{+\infty}\frac{c_0 q_0^{-\frac{1}{2}}\lambda_n^{-\frac{3}{2}}}{4\pi |x-z|}\,\left(\int_D e_n(y)\,dy\right)^2\,\int_0^{t-c_0^{-1}|x-z|} \sin\left[c_0(q_0\lambda_n)^{-\frac{1}{2}}(t-c_0^{-1}|x-z|-\tau)\right]\,V(z,\,\tau)\,d\tau \nonumber\\
&& \qquad\qquad\qquad\qquad\qquad\qquad + O(a^2) \nonumber\\
&=&\sum_{n=1}^{+\infty}\frac{c_0 q_0^{-\frac{1}{2}}\lambda_n^{-\frac{3}{2}}}{4\pi |x-z|}\,\left(\int_D e_n(y)\,dy\right)^2\,\int_{c_0^{-1}|x-z|}^{t} \sin\left[c_0(q_0\lambda_n)^{-\frac{1}{2}}(t-s)\right]\,V(z,\,s-c_0^{-1}|x-z|)\,ds \nonumber\\
&& \qquad\qquad\qquad\qquad\qquad\qquad + O(a^2)\nonumber\\
&=&\sum_{n=1}^{+\infty}\frac{c_1\lambda_n^{-\frac{3}{2}}}{4\pi |x-z|}\,\left(\int_D e_n(y)\,dy\right)^2\,\int_{c_0^{-1}|x-z|}^{t} \sin\left[ \frac{c_1}{\sqrt{\lambda_n}}(t-s)\right]\,V(z,\,s-c_0^{-1}|x-z|)\,ds \nonumber\\
&& \qquad\qquad\qquad\qquad\qquad\qquad + O(a^2), \quad (x,\,t) \in (\mathbb R^3\setminus\overline D)\times (0,\,T),
\end{eqnarray}
Thus, Theorem \ref{Asymptot-Introduction} is proved.

\section{Proof of Theorem \ref{Riesz-Introduction}}\label{sec-Riesz}

Based on the expansion \eqref{L-S_8}, we show how to reconstruct $V(z,\,t)$ from the measurement of $W(x,\,t)$ which is obtained by measuring the wave fields before and after injecting the particle $D$.  

Let $\{{\bar{\lambda}}_n\}_{n \in \mathbb{N}}$ be the sequence of eigenvalues of the Newtonian operator stated on the scaled particle $B$.
We start with the following comparison result shown in \cite{Suragan-2012}
\begin{equation}
{\bar{\lambda}}^N_n \leq {\bar{\lambda}}^{-1}_n \leq {\bar{\lambda}}^D_n,\; ~~ n\in \mathbb{N},
\label{Comparison}
\end{equation}
where $\{{\bar{\lambda}}^N_n\}_{n \in \mathbb{N}}$ and $\{{\bar{\lambda}}^D\}_{n \in \mathbb{N}}$ are the sequences of eigenvalues of the Laplacian stated on $B$ corresponding to Neumann and Dirichlet boundary conditions, respectively.

The following Weyl's expansions, in terms of the order $n$, of the eigenvalues are well known
\begin{equation}
{\bar{\lambda}}^N_n= \frac{6\pi^2}{\textrm{Vol}(B)}\; n^{2/3} +o(1) \; \mbox{ and } \; {\bar{\lambda}}^D_n= \frac{6\pi^2}{\textrm{Vol}(B)}\; n^{2/3} +o(1) \quad \mbox{ as } n\gg 1.
\end{equation}
Using (\ref{Comparison}), we deduce the similar properties for $\{{\bar{\lambda}}_n\}_{n \in \mathbb{N}}$
\begin{equation}
{\bar{\lambda}}^{-1}_n= \frac{6\pi^2}{\textrm{Vol}(B)}\; n^{2/3} +o(1) \quad \mbox{ as } n\gg 1.
\end{equation}
Hence, due to the scale $\lambda_n=\bar{\lambda}_n a^2$ and $c_1\sim a$, we deduce that
\begin{equation*}
\omega_n:= c_1 \lambda^{-\frac{1}{2}}= \sqrt{\frac{6\pi^2}{\textrm{Vol}(B)}}\; n^{1/3} +o(1).
\end{equation*}
However, we are only interested in eigenvalues $\bar \lambda_n$ for which the eigenfunctions $\bar e_n$ have nonzero average, i.e., $$\int_B \bar e_n(y)\,dy\not=0.$$ To look for those eigenvalues, we first consider the special case that the domain is a ball.

\begin{lemma}\label{eigen_N}
Let $\lambda_{lj}$ be the eigenvalues of the three-dimensional Newton potential in the ball $\{x:\,|x|<a\}$. Then we have
\begin{equation}\label{eigen_1}
\lambda_{l j}^{-1}=\frac{\left[\mu_{j}^{\left(l+\frac{1}{2}\right)}\right]^{2}}{a^{2}}, \quad l=0,\,1,\, \ldots, \; j=1,\,2,\, \ldots,
\end{equation}
where the $\mu_{j}^{\left(l+\frac{1}{2}\right)}$ are the roots of the transcendental equation
\begin{equation}\label{eigen_2}
(2 l+1) J_{l+\frac{1}{2}}\left(\mu_{j}^{\left(l+\frac{1}{2}\right)}\right)+\frac{\mu_{j}^{\left(l+\frac{1}{2}\right)}}{2}\left[J_{l-\frac{1}{2}}\left(\mu_{j}^{\left(l+\frac{1}{2}\right)}\right)-J_{l+\frac{3}{2}}\left(\mu_{j}^{\left(l+\frac{1}{2}\right)}\right)\right]=0.
\end{equation}
The eigenfunctions corresponding to each eigenvalue $\lambda_{l j},$ form a complete orthogonal system and can be represented in the form
\begin{equation}\label{eigen_3}
u_{ljm}=J_{l+\frac{1}{2}}\left(\sqrt{\lambda_{l j}^{-1}} r\right) Y_{l}^{m}(\varphi, \theta),
\end{equation}
where $J_{l+\frac{1}{2}}$ are the Bessel functions and
\begin{eqnarray}
&&Y_{l}^{m}(\varphi, \theta)=P_{l}^{m}(\cos \theta) \cos m \varphi, \quad m=0,1, \ldots, l, \label{eigen_31} \\
&&Y_{l}^{m}(\varphi, \theta)=P_{l}^{|m|}(\cos \theta) \sin |m| \varphi, \quad m=-1, \ldots,-l  \label{eigen_32}
\end{eqnarray}
for $l=0,1, \ldots$ are spherical functions. Here $P_{l}^{m}$ are the associated Legendre polynomials and $(r, \theta, \varphi)$ are the spherical coordinates with $\theta\in[0,\,\pi]$ and $\varphi\in[0,\,2\pi]$.
\end{lemma}

Since
\begin{equation*}
\int_{-1}^1 P_l^0(x)\,dx =0,\quad l\geq 1,
\end{equation*}
we can easily observe that the eigenfunctions \eqref{eigen_3} have nonzero average only for $m=0$ and $l=0$. The equation \eqref{eigen_2} with $l=0$ reduces to
\begin{equation}\label{eigen_4}
J_{\frac{1}{2}}\left(\mu_{j}^{\left(\frac{1}{2}\right)}\right)+\frac{\mu_{j}^{\left(\frac{1}{2}\right)}}{2}\left[J_{-\frac{1}{2}}\left(\mu_{j}^{\left(\frac{1}{2}\right)}\right)-
J_{\frac{3}{2}}\left(\mu_{j}^{\left(\frac{1}{2}\right)}\right)\right]=0.
\end{equation}
Note that
\begin{eqnarray}
&&J_{-\frac{1}{2}}(x) = \sqrt{\frac{2}{\pi x}}\,\cos x, \; J_{\frac{1}{2}}(x) = \sqrt{\frac{2}{\pi x}}\,\sin x, \quad x>0, \label{Bessel_1}\\
&&J_{\frac{3}{2}}(x) = \frac{1}{x} J_{\frac{1}{2}}(x) - J_{-\frac{1}{2}}(x) = \frac{1}{x} \sqrt{\frac{2}{\pi x}}\,\sin x -  \sqrt{\frac{2}{\pi x}}\,\cos x, \quad x>0. \label{Bessel_2}
\end{eqnarray}
Then we derive from \eqref{eigen_4} that
\begin{equation}\label{eigen_41}
\sin \left(\mu_{j}^{\left(\frac{1}{2}\right)}\right) + 2 \mu_{j}^{\left(\frac{1}{2}\right)}\, \cos \left(\mu_{j}^{\left(\frac{1}{2}\right)}\right) = 0,
\end{equation}
that is,
\begin{equation}\label{eigen_42}
\tan \left(\mu_{j}^{\left(\frac{1}{2}\right)}\right) = - 2 \mu_{j}^{\left(\frac{1}{2}\right)}.
\end{equation}
Note that $\tan x$ is monotonously increasing in $(j\pi-\pi/2,\,j\pi+\pi/2)$ and goes to $-\infty$ as $x\to j\pi-\pi/2$.
So the equation \eqref{eigen_42} has a unique root in each interval $(j\pi-\pi/2,\,j\pi+\pi/2)$ for $ j=1,\,2,\,\ldots$. Moreover, we can also easily see that
\begin{equation}\label{eigen_5}
\mu_{j}^{\left(\frac{1}{2}\right)}=j\pi-\frac{\pi}{2}+\gamma_j, \quad j=1,\,2,\,\ldots,
\end{equation}
where $\gamma_j$ monotonously decreases to zero as $j\to +\infty$. Hence, for the spherical domain $\{x:\,|x|<a\}$, the eigenvalues $\bar\lambda_n$ with $\int_B \bar e_n(y)dy\not=0$ have the following asymptotic property:
\begin{equation*}
\bar \lambda_n^{-1} = a^{-2}(n\pi-\frac{\pi}{2}+\gamma_n)^2,\quad n=1,\,2,\,\ldots.
\end{equation*}
Consequently, we have
\begin{equation*}
\omega_n := c_1 \lambda_n^{-\frac{1}{2}} = c_1 \frac{n\pi - \frac{\pi}{2} + \gamma_n}{a}, \quad n=1,\,2,\,\ldots.
\end{equation*}

For simplicity, we  let
\begin{equation}\label{Riesz_b}
b:=\frac{c_1\pi}{a}=1.
\end{equation}
Then we have
\begin{equation}\label{expansion_omega_n}
\omega_n=n-\frac{1}{2}+\frac{\gamma_n}{\pi},\quad n\in\mathbb N.
\end{equation}

In the following, we show that both $\{\cos(\omega_nt)\}_{n\in\mathbb N}$ and $\{\sin(\omega_nt)\}_{n\in\mathbb N}$ are Riesz basis in $L^2[0,\,\pi]$.  To this end, let us state the following results in \cite{Chr2016, X-V2001}.

\begin{lemma}\label{Riesz_cond}
Let $\mathcal H$ be a Hilbert space. A sequence $\{p_n\}_{n\in\mathbb Z}$ is a Riesz basis in $\mathcal H$ if and only if it is complete in $\mathcal H$ and there exist two positive constants $\mathcal A$ and $\mathcal B$ such that for every finite scalar sequence $\{c_n\}$
\begin{equation}\label{Riesz_ineq}
\mathcal A \sum_{n\in\mathbb Z} |c_n|^2 \leq \left\| \sum_{n\in\mathbb Z} c_n p_n \right\|^2 \leq \mathcal B \sum_{n\in\mathbb Z} |c_n|^2.
\end{equation}
\end{lemma}

\begin{lemma}\label{H-V_cos}
Let $\{\xi_n\},\,n\in\mathbb N_0:=\{0\}\cup\mathbb N$ be a sequence of nonnegative numbers with the property that $\xi_k\not=\xi_m$ for $k\not=m$ and of the form $\xi_n=n+\gamma+\gamma_n$ with $\gamma_n\in[-l,\,l]$ for sufficiently large $n$, where the constants $\gamma\in[0,\,\frac{1}{2}]$ and $l\in(0,\,\frac{1}{4})$ satisfy the condition
\begin{equation}\label{H-V_cond}
(1+\sin(2\gamma\pi))^{1/2}\,(1-\cos(l\pi)) + \sin(l\pi)<1.
\end{equation}
Then $\{\cos(\xi_nt)\}_{n\in\mathbb N_0}$ is a Riesz basis in $L^2[0,\,\pi]$.
\end{lemma}

\begin{lemma}\label{H-V_sin}
Let $\{\xi_n\},\,n\in\mathbb N$ be a sequence of positive numbers with the property that $\xi_k\not=\xi_m$ for $k\not=m$ and of the form $\xi_n=n-\gamma+\gamma_n$ with $\gamma_n\in[-l,\,l]$ for sufficiently large $n$, where the constants $\gamma\in[0,\,\frac{1}{2}]$ and $l\in(0,\,\frac{1}{4})$ satisfy \eqref{H-V_cond}.
Then $\{\sin(\xi_nt)\}_{n\in\mathbb N}$ is a Riesz basis in $L^2[0,\,\pi]$.
\end{lemma}

We are in a position to show the following result.

\begin{lemma}\label{Riesz_cs}
Both $\{\cos(\omega_nt)\}_{n\in\mathbb N}$ and $\{\sin(\omega_nt)\}_{n\in\mathbb N}$ are Riesz basis in $L^2[0,\,\pi]$.
\end{lemma}

{\bf Proof.}
First, we observe that $\omega_n=n-\frac{1}{2}+\frac{\gamma_n}{\pi},\, n\in\mathbb N$ can be rewritten as $\tilde \omega_n=n+\frac{1}{2}+\frac{\gamma_{n+1}}{\pi},\, n\in\mathbb N_0$. Second, we know from $\lim_{n\to+\infty}\gamma_n=0$ that $\frac{\gamma_n}{\pi}\in[-l,\,l]$ with $l\in(0,\,\frac{1}{4})$ for large enough $n$. In addition, we can easily check that the condition \eqref{H-V_cond} with $\gamma=1/2$ is satisfied. Hence, by Lemma \ref{H-V_cos}, we conclude that  $\{\cos(\omega_nt)\}_{n\in\mathbb N}$ is a Riesz basis in $L^2[0,\,\pi]$. Similarly, we can see from Lemma \ref{H-V_sin} that $\{\sin(\omega_nt)\}_{n\in\mathbb N}$ is also a Riesz basis in $L^2[0,\,\pi]$.
\hfill $\Box$

\medskip
Now define the sequence $\{p_n(t)\}_{n\in \mathbb Z^+\cup\mathbb Z^-}$ by
\begin{equation}\label{p_n}
p_n(t):=
\begin{cases}
\sin(\omega_nt), & n\in \mathbb Z^+, \\
\cos(\omega_{-n}t), & n\in\mathbb Z^-,
\end{cases}
\quad t\in[-\pi,\,\pi].
\end{equation}
Then we have the following key result.
\begin{theorem}\label{Riesz_pn}
The sequence $\{p_n(t)\}_{n\in \mathbb Z^+\cup\mathbb Z^-}$ is a Riesz basis in $L^2[-\pi,\,\pi]$.
\end{theorem}

{\bf Proof.} For any $f\in L^2[-\pi,\,\pi]$, we define
\begin{equation}\label{Riesz_sy_an}
f_{sy}(t):=\frac{f(t)+f(-t)}{2},\quad f_{an}(t):=\frac{f(t)-f(-t)}{2}.
\end{equation}
Note that $f_{sy}(t)$ is symmetric and $f_{an}(t)$ is antisymmetric with respect to $t$. Since $\{\cos(\omega_nt)\}_{n\in\mathbb N}$ is a Riesz basis in $L^2[0,\,\pi]$, we expand $f_{sy}(t)$ as
\begin{equation*}
f_{sy}(t) = \sum_{n=1}^{+\infty} a_n \cos(\omega_n t), \quad t\in [0,\,\pi].
\end{equation*}
Moreover, there exist two positive constants $\mathcal A_1$ and $\mathcal A_2$ such that
\begin{equation*}
\mathcal A_1 \sum_{n\in\mathbb Z^+} |a_n|^2 \leq \|f_{sy}\|^2_{L^2[0,\,\pi]} \leq \mathcal A_2 \sum_{n\in\mathbb Z^+} |a_n|^2.
\end{equation*}
By the symmetry of $f_{sy}(t)$ and $\cos(\omega_nt)$, we have
\begin{equation}\label{Riesz_p_1}
f_{sy}(t) = \sum_{n=1}^{+\infty} a_n \cos(\omega_n t), \quad t\in [-\pi,\,\pi]
\end{equation}
with
\begin{equation} \label{Riesz_p_12}
2\mathcal A_1 \sum_{n\in\mathbb Z^+} |a_n|^2 \leq \|f_{sy}\|^2_{L^2[-\pi,\,\pi]} = 2 \|f_{sy}\|^2_{L^2[0,\,\pi]} \leq 2\mathcal A_2 \sum_{n\in\mathbb Z^+} |a_n|^2.
\end{equation}
Similarly, using the antisymmetry of $f_{an}(t)$ and $\sin(\omega_nt)$ with respect to $t$, we can expand $f_{an}(t)$ as
\begin{equation}\label{Riesz_p_2}
f_{an}(t) = \sum_{n=1}^{+\infty} b_n \sin(\omega_n t), \quad t\in [-\pi,\,\pi]
\end{equation}
with
\begin{equation}\label{Riesz_p_22}
2\mathcal B_1 \sum_{n\in\mathbb Z^+} |b_n|^2 \leq \|f_{an}\|^2_{L^2[-\pi,\,\pi]} = 2 \|f_{an}\|^2_{L^2[0,\,\pi]} \leq 2\mathcal B_2 \sum_{n\in\mathbb Z^+} |b_n|^2,
\end{equation}
where $\mathcal B_1$ and $\mathcal B_2$ are two positive constants. 

From \eqref{Riesz_p_1} and \eqref{Riesz_p_2}, we obtain
\begin{equation}\label{Riesz_p_3}
f(t)=f_{sy}(t)+f_{an}(t) =  \sum_{n=1}^{+\infty} a_n \cos(\omega_n t) +  \sum_{n=1}^{+\infty} b_n \sin(\omega_n t) =:  \sum_{n\in\mathbb Z^+\cup\mathbb Z^-} c_n p_n(t), \quad t\in[-\pi,\,\pi],
\end{equation}
where
\begin{equation*}
c_n:=\begin{cases}
b_n, & n\in\mathbb Z^+,\\
a_{-n}, & n\in\mathbb Z^-.
\end{cases}
\end{equation*}
We observe from the definitions of $f_{sy}$ and $f_{an}$, i.e., \eqref{Riesz_sy_an}, that
\begin{equation*}
\|f_{sy}\|^2_{L^2[-\pi,\,\pi]}  \leq \|f\|^2_{L^2[-\pi,\,\pi]}, \quad \|f_{an}\|^2_{L^2[-\pi,\,\pi]}  \leq \|f\|^2_{L^2[-\pi,\,\pi]}.
\end{equation*}
In addition, we see from $f(t)=f_{sy}(t)+f_{an}(t) $ that
\begin{equation*}
\|f\|^2_{L^2[-\pi,\,\pi]}  \leq  2\left( \|f_{sy}\|^2_{L^2[-\pi,\,\pi]} + \|f_{an}\|^2_{L^2[-\pi,\,\pi]}\right).
\end{equation*}
Hence, we deduce from \eqref{Riesz_p_12} and \eqref{Riesz_p_22} that
\begin{eqnarray*}
\min\{\mathcal A_1,\,\mathcal B_1\} \sum_{n\in\mathbb Z^+\cup\mathbb Z^-} |c_n|^2
&=& \min\{\mathcal A_1,\,\mathcal B_1\}\left( \sum_{n\in\mathbb Z^-} |a_{-n}|^2 +  \sum_{n\in\mathbb Z^+} |b_n|^2\right) \\
&\leq& \mathcal A_1 \sum_{n\in\mathbb Z^+} |a_n|^2 + \mathcal B_1 \sum_{n\in\mathbb Z^+} |b_n|^2 \\
&\leq&\frac{1}{2}\left( \|f_{sy}\|^2_{L^2[-\pi,\,\pi]} +  \|f_{an}\|^2_{L^2[-\pi,\,\pi]}\right)\\
&\leq& \|f\|^2_{L^2[-\pi,\,\pi]}\\
&\leq&2\left( \|f_{sy}\|^2_{L^2[-\pi,\,\pi]} + \|f_{an}\|^2_{L^2[-\pi,\,\pi]}\right)\\
&\leq& 4\left(\mathcal A_2 \sum_{n\in\mathbb Z^+} |a_n|^2 +\mathcal  B_2 \sum_{n\in\mathbb Z^+} |b_n|^2\right) \\
&\leq& 4\max\{\mathcal A_2,\,\mathcal B_2\} \left( \sum_{n\in\mathbb Z^-} |a_{-n}|^2 + \sum_{n\in\mathbb Z^+} |b_n|^2 \right)\\
&=& 4\max\{\mathcal A_2,\,\mathcal B_2\} \sum_{n\in\mathbb Z^+\cup\mathbb Z^-} |c_n|^2.
\end{eqnarray*}
The proof is complete.  \hfill $\Box$

\bigskip
Let $\mathcal T_c,\,\mathcal T_s:\,L^2[0,\,\pi] \to L^2[0,\,\pi]$ be the operators corresponding to the Riesz basis $\{\cos(\omega_nt)\}_{n\in\mathbb N}$ and $\{\sin(\omega_nt)\}_{n\in\mathbb N}$, respectively. That is,
\begin{equation}\label{T_c}
\mathcal T_c[\cos((n-1)t)] = \cos(\omega_nt), \quad t\in[0,\,\pi],\; n\in\mathbb N
\end{equation}
and
\begin{equation}\label{T_s}
\mathcal T_s[\sin(nt)] = \sin(\omega_nt), \quad t\in[0,\,\pi],\; n\in\mathbb N.
\end{equation}
Define $L^2_{sy}[-\pi,\,\pi]$ and $L^2_{an}[-\pi,\,\pi]$ by
\begin{eqnarray*}
L^2_{sy}[-\pi,\,\pi]&:=&\left\{ f\in L^2[-\pi,\,\pi],\, f(-t) = f(t), \; a.e., \,t\in[-\pi,\,\pi]  \right\},\\
L^2_{an}[-\pi,\,\pi]&:=&\left\{ f\in L^2[-\pi,\,\pi],\, f(-t) = -f(t), \; a.e., \,t\in[-\pi,\,\pi]  \right\}.
\end{eqnarray*}
We now extend the operators $\mathcal T_c,\,\mathcal T_s$ to $L^2_{sy}[-\pi,\,\pi]$ and $L^2_{an}[-\pi,\,\pi]$. That is, for $f_{sy}\in L^2_{sy}[-\pi,\,\pi]$ and $f_{an}\in L^2_{an}[-\pi,\,\pi]$, we define
\begin{equation*}
\tilde {\mathcal T}_c [f_{sy}(t)] := 
\begin{cases}
\mathcal T_c [f_{sy}(t)] , & t\in[0,\,\pi],\\
\mathcal T_c [f_{sy}(-t)] , & t\in[-\pi,\,0]
\end{cases}
\end{equation*}
and
\begin{equation*}
\tilde {\mathcal T}_s [f_{an}(t)] := 
\begin{cases}
\mathcal T_s [f_{an}(t)] , & t\in[0,\,\pi],\\
-\mathcal T_s [f_{an}(-t)] , & t\in[-\pi,\,0],
\end{cases}
\end{equation*}
respectively.
Then we easily see that the operators
\begin{equation*}
\tilde {\mathcal T}_c:\, L^2_{sy}[-\pi,\,\pi] \to L^2_{sy}[-\pi,\,\pi], \quad \tilde {\mathcal T}_s:\, L^2_{an}[-\pi,\,\pi] \to L^2_{an}[-\pi,\,\pi]
\end{equation*}
are isomorphisms.

Finally, noticing $L^2[-\pi,\,\pi]=L^2_{sy}[-\pi,\,\pi]\textcircled+ L^2_{an}[-\pi,\,\pi]$,  we define the operator $\mathcal T:\,L^2[-\pi,\,\pi] \to L^2[-\pi,\,\pi]$ by $\mathcal T=\tilde {\mathcal T}_c \textcircled{+} \tilde {\mathcal T}_s$, which is also an isomorphism. Define
\begin{equation}\label{xi_n}
\xi_n(t):=
\begin{cases}
\sin(nt), & n\in\mathbb Z^+, \\
\cos((-n-1)t), &n\in\mathbb Z^-,
\end{cases}
\quad t\in [-\pi,\,\pi].
\end{equation}
Then we have
\begin{equation}\label{defn_T}
\mathcal T[\xi_n(t)] = p_n(t), \quad n\in\mathbb Z^+\cup\mathbb Z^-.
\end{equation}

%\bigskip

\section{Proof of Theorem \ref{Reconstruction-Introduction}}\label{sec-Reconstruction}
Based on the property \eqref{defn_T}, we reconstruct the field $V(z,\,t)$ from the measurement  $W(x, t)$ for a fixed $x$ outside $\Omega$ using the expansion (\ref{L-S_8}). The procedure is as follows.

We rewrite (\ref{L-S_8}) as
\begin{equation*}
W(x,\,t)=\; \sum_{n=1}^{+\infty}\frac{\omega_n \left(\int_D e_n(y)\,dy\right)^2}{4\pi |x-z| \lambda_n}\int_{c_0^{-1}|x-z|}^{t} \sin(\omega_n (t-s))\, V(z,\,s-c_0^{-1}|x-z|)\,ds + O(a^2),
\end{equation*}
or
\begin{eqnarray*}
W(x,\,t)&=& \sum_{n=1}^{+\infty}\frac{\omega_n \left(\int_D e_n(y)\,dy\right)^2}{4\pi |x-z| \lambda_n}\; \sin(\omega_n t) \int_{c_0^{-1}|x-z|}^{t} \cos(\omega_n\, s)\, V(z,\,s-c_0^{-1}|x-z|)\,ds \nonumber\\
&& -\sum_{n=1}^{+\infty}\frac{\omega_n \left(\int_D e_n(y)\,dy\right)^2}{4\pi |x-z| \lambda_n}\; \cos(\omega_n t) \int_{c_0^{-1}|x-z|}^{t} \sin(\omega_n\, s)\, V(z,\,s-c_0^{-1}|x-z|)\,ds + O(a^2).
\end{eqnarray*}
Note that the field $V$ can be expressed by a retarded potential
\begin{equation*}
V(x,\,t)=\dfrac{1}{4\pi c_0^2} \int_{B_x} \dfrac{J(y,\,t-|x-y|/c_0)}{|x-y|}\, dx,
\end{equation*}
where $B_x$ is a ball of radius $c_0 t$ centered at $x$. Due to the compactness of the support of $J$, we further suppose that the field $V$ is compactly supported with respect to t in $[0,\,\tilde T]$. Then, for $t>\tilde T+c_0^{-1} |x-z|$, the integrals in the above expansion of $W$ are independent of $t$. Define
\begin{equation*}
\alpha_n(x,\,z):= \frac{\omega_n \left(\int_D e_n(y)\,dy\right)^2}{4\pi |x-z| \lambda_n}.
\end{equation*}
Since the Riesz basis is defined on the interval $[-\pi, \pi]$, we need to shift the time interval. Set 
$$t_*:=\tilde T+c_0^{-1}|x-z|+\pi.$$ 
Using
\begin{eqnarray*}
\sin(\omega_nt) &=& \sin[\omega_n(t-t_*)+\omega_nt_*] = \sin(\omega_n(t-t_*)) \, \cos(\omega_nt_*) + \cos(\omega_n(t-t_*)) \, \sin(\omega_nt_*), \\
\cos(\omega_nt) &=& \cos[\omega_n(t-t_*)+\omega_nt_*] = \cos(\omega_n(t-t_*)) \, \cos(\omega_nt_*)  - \sin(\omega_n(t-t_*)) \, \sin(\omega_nt_*),
\end{eqnarray*}
we have
\begin{eqnarray}\label{Riesz_h1}
W(x,\,t)  &=& \sum_{n=1}^{+\infty}\alpha_n(x,\,z)\, \sin(\omega_n (t-t_*)) \int_{c_0^{-1}|x-z|}^{\tilde T +c_0^{-1}|x-z|} \cos(\omega_nt_*)\,\cos(\omega_n s)\, V(z,\,s-c_0^{-1}|x-z|)\,ds \nonumber\\
&& +  \sum_{n=1}^{+\infty}\alpha_n(x,\,z)\, \cos(\omega_n (t-t_*)) \int_{c_0^{-1}|x-z|}^{\tilde T +c_0^{-1}|x-z|} \sin(\omega_nt_*)\,\cos(\omega_n s)\, V(z,\,s-c_0^{-1}|x-z|)\,ds \nonumber\\
&& -\sum_{n=1}^{+\infty}\alpha_n(x,\,z)\, \cos(\omega_n (t-t_*)) \int_{c_0^{-1}|x-z|}^{\tilde T +c_0^{-1}|x-z|} \cos(\omega_nt_*)\sin(\omega_n s)\, V(z,\,s-c_0^{-1}|x-z|)\,ds  \nonumber\\
&& +\sum_{n=1}^{+\infty}\alpha_n(x,\,z)\, \sin(\omega_n (t-t_*)) \int_{c_0^{-1}|x-z|}^{\tilde T +c_0^{-1}|x-z|} \sin(\omega_nt_*) \sin(\omega_n s)\, V(z,\,s-c_0^{-1}|x-z|)\,ds \nonumber \\
&& + O(a^2) \nonumber\\
&=& \sum_{n=1}^{+\infty}\alpha_n(x,\,z)\, \sin(\omega_n (t-t_*)) \int_{c_0^{-1}|x-z|}^{\tilde T +c_0^{-1}|x-z|} \cos(\omega_n(s-t_*))\,V(z,\,s-c_0^{-1}|x-z|)ds \nonumber\\
&& -\sum_{n=1}^{+\infty}\alpha_n(x,\,z)\, \cos(\omega_n (t-t_*)) \int_{c_0^{-1}|x-z|}^{\tilde T +c_0^{-1}|x-z|} \sin(\omega_n\, (s-t_*))\, V(z,\,s-c_0^{-1}|x-z|)\,ds \nonumber\\
&& + O(a^2).
\end{eqnarray}
Now we measure $W(x,\,t)$ for $t\in[\tilde T+c_0^{-1}|x-z|,\,\tilde T+c_0^{-1}|x-z|+2\pi]$. Set $\tau:=t-t_*.$
Applying the operator $\mathcal T^{-1}$ to \eqref{Riesz_h1}, we have
\begin{eqnarray}\label{Riesz_h2}
&&\mathcal T^{-1}[W(x,\,t_*+\tau)] \nonumber\\
&=&  \sum_{n=1}^{+\infty}\alpha_n(x,\,z)\, \sin(n \tau) \int_{c_0^{-1}|x-z|}^{\tilde T +c_0^{-1}|x-z|} \cos(\omega_n(s-t_*))\,V(z,\,s-c_0^{-1}|x-z|)\,ds \nonumber\\
&& -\sum_{n=1}^{+\infty}\alpha_n(x,\,z)\, \cos((n-1) \tau ) \int_{c_0^{-1}|x-z|}^{\tilde T +c_0^{-1}|x-z|} \sin(\omega_n\, (s-t_*))\, V(z,\,s-c_0^{-1}|x-z|)\,ds \nonumber\\
&& + O(a^2).
\end{eqnarray}
As the sequence $\{\cos((n-1)t),\,\sin(nt) \}_{n=1}^\infty$ is the Fourier basis in $L^2[-\pi,\,\pi]$, we  can approximately get the coefficients
\begin{equation}\label{Riesz_h3}
A_n(x,\,z):=\int_{c_0^{-1}|x-z|}^{\tilde T +c_0^{-1}|x-z|} \cos(\omega_n (s-t_*))\,V(z,\,s-c_0^{-1}|x-z|)\,ds, \quad n\in\mathbb N
\end{equation}
and
\begin{equation}\label{Riesz_h4}
B_n(x,\,z):=\int_{c_0^{-1}|x-z|}^{\tilde T +c_0^{-1}|x-z|} \sin(\omega_n\, (s-t_*)) \,V(z,\,s-c_0^{-1}|x-z|)\,ds, \quad n\in\mathbb N.
\end{equation}
Suppose that 
\begin{equation}\label{assume_T}
0< \tilde T \leq 2\pi.
\end{equation}
Set $\tau_*:=s-t_*+\tilde T.$ Then we further have
\begin{eqnarray}\label{Riesz_h5}
A_n(x,\,z)&=&\int_{c_0^{-1}|x-z|}^{2\pi +c_0^{-1}|x-z|} \cos(\omega_n (s-t_*))\,V(z,\,s-c_0^{-1}|x-z|)\,ds \nonumber\\
&=&\int_{-\pi}^{\pi} \cos(\omega_n (\tau_*- \tilde T)) \,V(z,\,\tau_*+\pi)\,d\tau_* \nonumber\\
&=&\cos(\omega_n\tilde T ) \int_{-\pi}^{\pi} \cos(\omega_n \tau_*) \,V(z,\,\tau_*+\pi)\,d\tau_* \nonumber\\
&& + \sin(\omega_n\tilde T ) \int_{-\pi}^{\pi} \sin(\omega_n \tau_*) \,V(z,\,\tau_*+\pi)\,d\tau_*, \quad n\in\mathbb N
\end{eqnarray}
and
\begin{eqnarray}\label{Riesz_h6}
B_n(x,\,z)&:=&\int_{c_0^{-1}|x-z|}^{2\pi +c_0^{-1}|x-z|} \sin(\omega_n\, (s-t_*)) \,V(z,\,s-c_0^{-1}|x-z|)\,ds \nonumber\\
&=&\int_{-\pi}^{\pi} \sin(\omega_n (\tau_*- \tilde T)) \,V(z,\,\tau_*+\pi)\,d\tau_* \nonumber\\
&=&\cos(\omega_n\tilde T ) \int_{-\pi}^{\pi} \sin(\omega_n \tau_*) \,V(z,\,\tau_*+\pi)\,d\tau_* \nonumber\\
&& - \sin(\omega_n\tilde T ) \int_{-\pi}^{\pi} \cos(\omega_n \tau_*) \,V(z,\,\tau_*+\pi)\,d\tau_*, \quad n\in\mathbb N.
\end{eqnarray}
Hence, we can easily obtain the following two quantities:
\begin{eqnarray*}
C_n(x,\,z)&:=& \int_{-\pi}^{\pi} \cos(\omega_n \tau_*) \,V(z,\,\tau_*+\pi)\,d\tau_*,\quad n\in\mathbb N, \\
D_n(x,\,z)&:=& \int_{-\pi}^{\pi} \sin(\omega_n \tau_*) \,V(z,\,\tau_*+\pi)\,d\tau_*,\quad n\in\mathbb N,
\end{eqnarray*}
since \eqref{Riesz_h5} and \eqref{Riesz_h6} constitute a linear system with respect to $C_n(x,\,z)$ and $D_n(x,\,z)$. By Corollary 3.6.3 in \cite{Chr2016}, we get
\begin{equation}\label{Riesz_h70}
V(z,\,t+\pi) = \sum_{n=1}^{+\infty}\left[C_n(x,\,z)g_n(t) + D_n(x,\,z) h_n(t)\right], \quad t\in[-\pi,\,\pi],
\end{equation}
where $g_n(t):=(\mathcal T^{-1})^*\cos((n-1) t)$ and  $h_n(t):=(\mathcal T^{-1})^*\sin(n t)$.

\medskip
By varying the location of $z$ in $\Omega$, we can get the field $V(\cdot,\,t)|_\Omega$ for $t\in[0,\,2\pi]$. Observe that as we measure $W(x, t)$ for the time $t$ larger than $T_J$, we have $V(x,\, t)=0$, and hence we actually measure only $U(x,\, t)$. Then, we can reconstruct the source $J(x,\,t)$ from the measurement data of $U(x,\,t)$ for a fixed $x$ over a time interval, by numerical differentiation applied to $V(\cdot,\,t)|_\Omega$, for $t\in[0,\,2\pi]$, and the equation $c_0^{-2} V_{tt} - \Delta V = J$.

Here, we would like to emphasize that the condition \eqref{assume_T} is not critical, and the argument can work for any fixed $\tilde T$. Indeed, in the above derivations, we simply let $b= 1$ such that we could work on the standard interval $[-\pi,\, \pi]$. We can easily see that $2\pi$ is just the length of the time interval, nothing else. For general $b$, we can construct the Riesz basis on the interval $[-\pi/b, \pi/b]$ and do the shift similarly.

\bigskip
{\bf Acknowledgement:} This work is supported by National Natural Science Foundation of China (Nos. 12071072, 11971104) and the Austrian Science Fund (FWF): P 30756-NBL.

\end{document}